\documentclass[12pt]{amsart}
\usepackage{a4}
\usepackage[cp1250]{inputenc}
\usepackage{amssymb}
\usepackage[T1]{fontenc}

\newcounter{licz}

\newenvironment{wylicz}{\begin{list}%
{{\rm(\alph{licz})\hfill}}{%
\usecounter{licz}%
\addtolength{\leftmargin}{-.45\leftmargin}%
\addtolength{\labelwidth}{\leftmargin}%
\addtolength{\labelwidth}{-1\labelsep}%
\addtolength{\topsep}{1.5\topsep}%
}}{\end{list}}

\newenvironment{lista}{\begin{list}%
{{\rm---\hfill}}{%
\setlength{\labelsep}{0ex}%
\settowidth{\labelwidth}{---}%
\setlength{\leftmargin}{\parindent}%
\addtolength{\topsep}{1.5\topsep}%
}}{\end{list}}

\newcounter{stat}[section]

\newtheorem{statem}{}[section]

\newcommand{\h}{\hspace*{-.08ex}}
\newcommand{\g}{\hspace*{-.15ex}}

\newcommand{\id}{\mathop{\mathrm{id}}\nolimits}
\newcommand{\diam}{\mathop{\mathrm{diam}}\nolimits}
\newcommand{\fr}{\mathop{\mathrm{bd}}\nolimits}
\newcommand{\cl}{\mathop{\mathrm{cl}}\nolimits}

\newcommand{\card}{\mathop{\mathrm{card}}\nolimits}
\renewcommand{\dim}{\mathop{\mathrm{dim}}\nolimits}
\newcommand{\Ind}{\mathop{\mathrm{Ind}}\nolimits}
\newcommand{\Indo}{\mathop{\mathrm{Ind}_0}\nolimits}
\newcommand{\ind}{\mathop{\mathrm{ind}}\nolimits}

\newcommand{\ball}[2]{\mathop{\mathrm{B\mbox{}}}\nolimits (#1,#2)}

\title[Fully closed maps and Anderson-Choquet continua]{Fully closed
maps and non-metrizable higher-dimensional Anderson-Choquet continua}

\author[J.\ Krzempek]{Jerzy Krzempek}

\thanks{The author was partially supported by MNiSW grant no.\ N201 034 31/2717}

\address{Institute of Mathematics\\
Silesian University of Technology\\
Kaszubska 23\\
44-100 Gliwice\\
Poland
\newline\indent
(Instytut Matematyki\\
Politechnika Śląska)}

\email{j.krzempek@polsl.pl}



\begin{document}

\begin{abstract}
Fedorchuk's fully closed  (continuous) maps and resolutions are
applied in constructions of non-metrizable higher-dimensional
analogues of Anderson, Choquet, and Cook's rigid continua. Certain
theorems on dimension-lowering maps are proved for inductive
dimensions and fully closed maps from spaces that need not be
hereditarily normal, and some examples of continua have
non-coinciding dimensions.
\end{abstract}

\maketitle

Fully closed (continuous) maps and resolutions appear in numerous constructions (see S. Watson \cite{wat}, V.V.\ Fedorchuk
\cite{fed} for surveys), in particular, in constructions of homogeneous spaces with non-coinciding dimensions. In this paper
we apply the maps in order to obtain examples of continua with strong, hereditary\linebreak rigidity properties.

A non-degenerate continuum $X$ is called
\begin{lista}
\item[---]{\em an Anderson-Choquet continuum}\/ if every non-degenerate subcontinuum
$P$ of~$X$ has exactly one embedding $P\to X$, the identity
$\id_P$;
\item[---]{\em a Cook continuum}\/
if every non-degenerate subcontinuum $P$ of $X$ has exactly one
non-constant map $P\to X$, the identity $\id_P$.
\end{lista}
Examples were constructed---respectively---by R.D.\ Anderson, G.\
Choquet \cite{and} (a plane hereditarily decomposable continuum),
and H.\ Cook \cite{cook} (a metric, one-dimensional, hereditarily
indecomposable continuum). T.\ Maćkowiak \cite{mac15} constructed
a metric, chainable, hereditarily decomposable Cook continuum.

All known examples of such continua are one-dimensional. A metric
Cook continuum must have dimension $\leq 2$ (Mać\-ko\-wiak
\cite{mac2}), and if it is hereditarily indecomposable, then it
must be one-dimensional (Krzempek \cite{krz1}). On the other hand,
several authors investigated rigidity properties of
higher-dimensional continua (J.J.\ Chara\-to\-nik \cite{charat},
M.\ Reńska \cite{ren}, E.\ Pol [\ref{po1}--\ref{po4}, \ref{kep}],
see \cite{krz1} for more references). In \cite{krz1} the present
author constructed a metric, $n$-dimensional (arbitrary $n>1$),
hereditarily indecomposable continuum no two of whose disjoint
$n$-dimensional sub\-con\-tin\-ua are homeomorphic, but he was not
able to ensure that the continuum be Anderson-Choquet.

In this paper we achieve better results for the non-metric case:
we construct examples of non-metrizable higher-dimensional (with
respect to the dimensions dim, ind, and Ind) Anderson-Choquet
continua and Cook continua. Our other aim is a study of behavior
of the Charalambous-Filippov-Ivanov inductive dimension $\Indo$
under fully closed maps.

In Sections 1--3 we gather some facts about fully closed maps,
ring-like maps, and dimensions. Application of $\Indo$ simplifies
estimating inductive dimensions in some well-known examples. We
show that, {\em if $f$ is a fully closed map from a non-empty
compact space $X$ to a first countable space, then\/ $\Indo
X\leq\Indo fX + \Indo f$.} This enables us to prove that (1)
Fedorchuk's first countable compact spaces \cite{fed5} with
$\dim\!=\!n\!<\!2n-1\!\leq\!\ind\!\leq\! 2n$ have also
$\Ind\!\leq\!\Indo\!=\!2n$, and (2)~V.A.~Chatyrko's chainable
continua and homogeneous continua \cite{chat5} with $\dim= 1$ and
$\ind=n$ have also $\Ind =\Indo=n$. In Section~3 we slightly
modify the Fedorchuk-Emeryk-Chatyrko resolution theorem, and that
is our main tool for constructions.

Section 4 contains a construction of hereditarily indecomposable Anderson-Cho\-quet continua with $\dim =n$ (arbitrary
$n>1$) and a construction of a Cook continuum with $\dim =2$. In Section 5 we obtain chainable (hence, $\dim=1$),
hereditarily decomposable Cook continua with $n\leq\ind \leq\Ind\leq\Indo=n+1$. All the continua are separable and first
countable, and some have $\dim<\ind$.

\section{Preliminaries: continua, maps, and covering dimension}

By a {\em space} we mean a regular $T_1$ topological space, and all considered {\em maps} are continuous and closed. A {\em
continuum} is a non-empty connected compact space. A subcontinuum $A$ of a space $X$ is said to be {\em terminal}\/ if for
every continuum $B\subset X$, either $A\cap B=\emptyset$, $A\subset B$, or $B\subset A$. A continuum is said to be
\begin{lista}
\item[---]{\em decomposable} if it is the union of two proper subcontinua;
\item[---]{\em hereditarily decomposable}
(abbrev.\ HD) if each of its non-degenerate subcontinua is
decomposable;
\item[---]{\em hereditarily indecomposable}
(abbrev.\ HI) if none of its subcontinua is decomposable
(equivalently: each of its subcontinua is terminal); and
\item[---]{\em chainable} if for every open cover, there exists a natural number $n$ and a
closed refinement $F_1,F_2,\ldots,F_n$ of the cover such that
$F_i\cap F_j\neq\emptyset$ iff $|i-j|\leq 1$.
\end{lista}

A (closed continuous) map $f\colon X\to Y$ is said to be
\begin{lista}
\item[---]{\em irreducible}\/ if for every closed proper subset
$F\varsubsetneq X$, also $fF\varsubsetneq fX$;
\item[---]{\em monotone}\/ [\hspace{-.15ex}{\em
atomic}\hspace{.1ex}] if for every point $y\in Y$, the pre-image $f^{-1}y$ is a continuum [respectively: a terminal
continuum];
\item[---]{\em ring-like} if for every point $x\in X$ and
every pair of open sets $U\ni x$ and $V\ni fx$, there is an open set $V'\ni fx$ such that $V'\subset V$ and $f^{-1}\fr
V'\subset U$;
\item[---]{\em fully closed}\/ if for every pair of disjoint closed subsets
$F,G\subset X$, the intersection $fF\cap fG$ is a discrete subspace of $Y$.\footnote{An extensive survey \cite{fed} by
Fedorchuk is devoted to fully closed maps, ring-like maps, and their applications. See \cite[Section II.1]{fed} for
equivalent definitions of fully closed maps. For~terms not explicitly defined herein the reader is directed to the
monographs by R.\ Engelking \cite{eng2,eng1} and K.\ Kuratowski \cite{kur}.}
\end{lista}

We shall frequently use the simple fact that, {\em if $f$ is a fully closed \mbox{{\em
[}\hspace{-.07ex}ring-like\hspace{.15ex}{\em ]}} map from a space\/ $X$ and\/ $X'\subset X$ is closed, then also the
restriction $f|X'$ is fully closed\/ {\em [\hspace{-.1ex}}respectively: ring-like}\hspace{.1ex}].

The following proposition is well-known (cf.\ \cite[Proposition II.3.10]{fed}).

\begin{statem}\label{met}
{\bf Proposition.} Suppose that\/ $X$ and\/ $Y$ are compact
spaces, and\/ $f\colon X\to Y$ is a map whose every point-inverse
is metrizable.

If the set\/ $C_2f=\{y\in Y\colon\card f^{-1}y>1\}$ is countable
and\/ $Y$ is metrizable, then $X$ is metrizable. The converse is
true if moreover $f$ is a fully closed map.
\end{statem}

\begin{proof}
\g If \h $C_2f$ \h is \h countable \h and \h $Y$ \h is metrizable, \h then \h the \h diagonal \h$\{(x,x)\!\in\! X\!\times\!
X$: $x\in X\}$ is a $G_\delta$-subset of $X\times X$. A Shne\u\i der theorem (see R.\ Engelking \cite[Exercise 4.2.B]{eng2})
implies that $X$ is metrizable.

Assume that $f$ is fully closed and there is a metric on $X$.
Using sequential compactness of $X$, one easily checks that the
set $\{y\in Y\colon\diam f^{-1}y\geq 1/n\}$ is finite for every
$n$. Thus, $C_2f$ is countable.
\end{proof}

\begin{statem}\label{preat} {\bf Proposition} \em (A.\ Emeryk,
Z.\ Horbanowicz \cite[Theorem 1]{eh}). \em A~map $f$ from a continuum $X$ is atomic iff $A=f^{-1}fA$ for every continuum
$A\subset X$ such that $fA$ is not a single point.\qed
\end{statem}

The point in the foregoing proposition is that, if the
irreducibility condition $A=f^{-1}fA$ is satisfied for all
subcontinua $A\subset X$ with non-degenerate images, then $f$ is a
monotone map.

\begin{statem}\label{irr}\em
It is easily seen that any ring-like map $f\colon X\to Y$ has even a stronger property: it is {\em connected irreducible,}
i.e.\ $A=f^{-1}fA$ for every closed subspace $A\subset X$ such that $fA$ is connected and contains more than one point (see
\cite[II.1.15]{fed}).\linebreak In particular, $f$ is irreducible whenever $Y=fX$ is connected and contains more than one
point.
\end{statem}

\begin{statem}\label{at}
{\bf Proposition.} If $f$ is a ring-like map from a compact space
$X$ onto a non-degenerate continuum, then $X$ is a continuum and
$f$ is an atomic map.
\end{statem}

\begin{proof}
Assume that $f\colon X\to Y$ is a ring-like map, $X$ is compact,
and $fX$ is a non-degenerate continuum. Suppose that $X=A\cup B$,
where $A$ and $B$ are non-empty, closed, and disjoint. By the
remark \ref{irr}, the complement $fX\setminus fA$ is non-empty.
Consider any component $M$ of $fX\setminus fA$. By Janiszewski's
boundary bumping theorem (see K.\ Kuratowski \cite[\S 47, III,
Theorem 2]{kur}), the closure $\cl{M}$ meets $fA$, and is not a
single point. Since $\cl{M}\subset fB$ and $f$ is connected
irreducible, we obtain $f^{-1}\cl{M}\subset B$. A contradiction.
Therefore, $X$~is a continuum. Finally, $f$~is an atomic map by
Proposition \ref{preat} and the remark \ref{irr}.
\end{proof}

The next useful proposition belongs to folklore (see Fedorchuk \cite[the proof of Lemma III.3.3]{fed} and \cite[Proposition
1.3]{fed02}).

\begin{statem}\label{folk}
{\bf Proposition.} Suppose that $f$ is a \mbox{{\em
(}\hspace{-.2ex}closed{\em\hspace{.35ex})}} monotone map from a
space~$X$, and\/ $A,B\!\subset\! fX$ are disjoint closed sets. If
$L$ is a partition\/\footnote{We say that a closed set $L\subset
X$ is a {\em partition in $X$}\/ if the complement $X\setminus L$
is not connected. Moreover, $L$ is a {\em partition in $X$ between
disjoint sets $A,B\subset X$} if there are disjoint open sets
$U,V\subset X$ such that $X\setminus L=U\cup V$, $A\subset U$, and
$B\subset V$.} in $X$ between $f^{-1}A$ and $f^{-1}B$, then $fL$
is a partition in\/ $fX$ between $A$ and~$B$.\qed
\end{statem}

\begin{statem}\label{comp}
{\bf Proposition.} Suppose that $f\colon X\to Y$ and $g\colon Y\to
Z$ are surjective ring-like maps between compact spaces $X$, $Y$,
and~$Z$. If $g$ is a monotone map, then the composition $gf$ is
ring-like.
\end{statem}

\begin{proof}
Take $x\!\in\! X$ and open sets $U\!\ni\! x$, $W\!\ni\!
z\!=\!gfx$. We can assume $gfU\!\subset\! W$. Since $f$ is
ring-like, there is an open set $V'\ni fx$ such that $\cl
V'\subset g^{-1}W$ and $f^{-1}\fr V'\subset U$. If $g^{-1}z$ is
not a singleton, we can moreover have $g^{-1}z\not\subset\cl V'$.
Clearly, $\fr V'\subset Y\setminus f(X\setminus U)$. There are two
cases. (1) If $g^{-1}z$ is a non-degenerate continuum, then there
is a point $y\in\fr V'$ with $gy=z$. The set $Y\setminus
f(X\setminus U)$ is an open neighborhood of $y$, and there is an
open set $W'\ni z$ such that $W'\subset W$ and $g^{-1}\fr
W'\subset Y\setminus f(X\setminus U)$. Thus, $(gf)^{-1}\fr
W'\subset U$. (2)~If~$g^{-1}z$ is a single point, then $\fr V'$ is
a partition in $Y$ between $fx$ and $g^{-1}(Z\setminus W)$. By
Proposition \ref{folk}, $g\fr V'$ is a partition in $Z$ between
$z$ and $Z\setminus W$. Hence, there is an open set $T\subset\cl T
\subset W$ such that $z\in T$ and $\fr T\subset g\fr V'$. As $g$
is ring-like and $\fr T\subset W\cap g[Y\setminus f(X\setminus
U)]$, for each $t\in\fr T$ there is an open set $W_t\subset W$
such that $g^{-1}\fr W_t\subset Y\setminus f (X\setminus U)$.
Then, $\fr T\subset W_{t_1}\cup\ldots\cup W_{t_n}$, where
$t_1,\ldots,t_n\in\fr T$. We put $W'=T\cup W_{t_1}\cup\ldots\cup
W_{t_n}\subset W$, and have $\fr W'\subset \fr
W_{t_1}\cup\ldots\cup \fr W_{t_n}$.\linebreak Thus, we obtain
$g^{-1}\fr W'\subset Y\setminus f(X\setminus U)$ and $(gf)^{-1}\fr
W'\subset U$.
\end{proof}

\begin{statem}\label{partit}
{\bf Proposition} \em (cf.\ Chatyrko \cite[Proposition 2]{chat5}).
\em Suppose that $f\colon X\to Y$ is a surjective ring-like map
from a compact space $X$, and $g\colon Y\to Z$ is a surjective
monotone map without degenerate point-inverses. If every partition
in\/ $Y$ contains a point-inverse of\/ $g$, then every partition
in\/ $X$ contains a point-inverse of the composition\/ $gf$.
\end{statem}

\begin{proof}
Assume that $X\neq\emptyset$ and every partition in $Y$ contains a point-inverse of~$g$. Since the empty set is not a
partition in $Y$, $Y$~is a non-degenerate continuum. Hence, $f$ is irreducible by \ref{irr}, and monotone by Proposition
\ref{at}. Take a partition $L$ in $X$. The irreducibility of $f$ implies that $L$ is a partition between some point-inverses
$f^{-1}a$ and $f^{-1}b$, where $a,b\in Y$. By Proposition \ref{folk}, $fL$ is a partition in $Y$ between $a$ and $b$. Then
$g^{-1}z\subset fL$ for some $z\in Z$, and again by \ref{irr}, $f^{-1}g^{-1}z\subset L$.
\end{proof}

\begin{statem}\label{ch1}
{\bf Proposition}  \em (\hspace{-.2ex}\em implicite\/ \em Chatyrko
\cite{chat1}). \em If $f$ is a ring-like map from a compact space
$X$, then\/ $\dim fX\leq\dim X$.
\end{statem}

\begin{proof}[Proof \rm (cf.\ Chatyrko \cite{chat1}, p.\ 124).]
We shall prove that for every natural number~$n$, the inequality
$n\leq\dim fX$ implies $n\leq\dim X$. For $n=0$, this is obvious.
For $n=1$, $fX$ contains a non-degenerate continuum $Y$. Then, by
Proposition \ref{at}, $f^{-1}Y$ is a non-degenerate continuum, and
hence, $1\leq\dim X$.

Let us recall that a normal space $Y$ has $\dim Y\geq n$ iff there
exists an essential family $(A_1,B_1),\ldots,(A_n,B_n)$ in $Y$,
i.e.\ $A_i, B_i\subset Y$ are disjoint closed subsets for each
$i$, and for every partitions $L_i$ between $A_i$ and $B_i$, the
intersection $\bigcap_{i=1}^n L_i$ is non-empty (cf.\
En\-gel\-king \cite[Theorem 3.2.6]{eng1}).

Let $2\leq n\leq\dim fX$. Since $fX$ contains a component of dimension $\geq n$, we can assume that $fX$ is a continuum.
Then, $f$ is a monotone map by Proposition \ref{at}. Take an essential family $(A_1,B_1),\ldots,(A_n,B_n)$ in $fX$. We shall
show that the pre-images $(f^{-1}A_1,f^{-1}B_1),\ldots,(f^{-1}A_n,f^{-1}B_n)$ form an essential family in~$X$. If $L_i$ are
partitions in $X$ between $f^{-1}A_i$ and $f^{-1}B_i$, then $fL_i$ are partitions in $fX$ between $A_i$ and $B_i$
(Proposition \ref{folk}). By Lemma 5.2 in \cite{rsw}, the intersection $\bigcap_{i=2}^nfL_i$ contains a continuum $P$ which
meets both $A_1$ and $B_1$. Since $f(L_i\cap f^{-1}P)=P$ for $i=2,\ldots,n$, the remark \ref{irr} implies that
$f^{-1}P=L_i\cap f^{-1}P$ and $f^{-1}P\subset \bigcap_{i=2}^n L_i$. As $f$~is monotone, $f^{-1}P$ is a continuum, $f^{-1}P$
meets~$L_1$, and hence, $\bigcap_{i=1}^n L_i$ is non-empty. Therefore, $n\leq \dim X$.
\end{proof}

The fiberwise covering dimension of a map $f\colon X\to Y$ is
defined as
\[\dim f=\sup\{\dim f^{-1}y\colon y\in Y\}.\]
Other fiberwise dimension functions $\ind$, $\Ind$, etc for maps
are defined similarly.

\begin{statem}\label{fedo}
{\bf Theorem }\em (
Fedorchuk, see \cite[Theorem III.2.4]{fed}). \em If $f$ is a fully closed map from a normal space~$X$, then $\dim
X\leq\max\{\dim fX,\dim f\}$.\qed
\end{statem}

The following is a consequence of Theorem \ref{fedo} and
Proposition \ref{ch1}.

\begin{statem}\label{rdim}
{\bf Corollary.} If $f$ is a ring-like fully closed map from a
compact space $X$, then\/ $\dim X=\max\{\dim fX,\dim f\}$.\qed
\end{statem}

\section{Maps that reduce inductive dimensions}

Since the theory of Ind is unsatisfactory outside the class of
hereditarily normal spaces, we shall use another inductive
dimension function $\Indo$, which was introduced by M.G.\
Charalambous \cite{char1,char2} and A.V.\ Ivanov \cite{iv}.

\begin{statem} {\bf Definition.} \em
For normal spaces $X$, $\Indo X\in\{-1,0,1,2,\ldots,\infty\}$ is
defined so that
\begin{wylicz}
\item $\Indo X=-1$ iff $X=\emptyset$;
\item $\Indo X\leq n\geq 0$ iff for every pair of disjoint closed
sets $A,B\subset X$, between $A$ and~$B$ there is a $G_\delta$
partition $L$ such that $\Indo L\leq n-1$;
\item $\Indo X=n$ iff $\Indo X\leq n$ and it is not true that
$\Indo X\leq n-1$;
\item $\Indo X=\infty$ if for every $n\in \mathbb{N}$, it is not
true that $\Indo X\leq n$.
\end{wylicz}
\end{statem}

It is clear that $\Ind X\leq\Indo X$ for every normal space $X$,
and $\Ind X=\Indo X$ if $X$~is perfectly normal.

\begin{statem}\label{count}\em
{\bf Countable sum theorem for} $\Indo$ (Charalambous \cite{char2}, Ivanov \cite{iv}). \em Suppose that\/
$X=\bigcup_{i=1}^\infty F_i$ is a normal space, and $F_i$ are closed $G_\delta$-subsets of $X$. If\/ $\Indo F_i\leq n$ for
every\/ $i$, then $\Indo X\leq n$.\qed
\end{statem}

The assumption that $F_i$ are $G_\delta$-sets is necessary in Theorem \ref{count} even if $X$ is a hereditarily normal
compact space, see \cite{iv}. Besides \cite{char2,iv}, see Charalambous and Cha\-tyr\-ko \cite{chacha} for more (also
bibliographical) information on $\Indo$.

The following theorem on dimension-lowering fully closed maps
seems to be important because of its applications.

\begin{statem}\label{indor}
{\bf Theorem.} If $f$ is a fully closed map from a non-empty
normal space $X$ to a space whose every discrete closed subspace
is a $G_\delta$-set, then $\Indo X\leq\Indo fX + \Indo f$.
\end{statem}

We shall modify the proof of Theorem III.2.8 in \cite{fed}. At
first, we need some standard preparation (see \cite[pp.\
4213--4216]{fed} for details). Let $f\colon X\to Y$ be a map, and
$M\subset Y$ be an arbitrary set. Consider the decomposition
$$\mathcal{M}=\{f^{-1}y\colon y\in Y\setminus M\}\cup \{\{x\}\colon x\in f^{-1}M\}$$
of $X$. Let $Y^M=X/\mathcal{M}$ be the quotient space, $f^M\colon
X\to Y^M$ the natural quotient projection, and $\pi^M\colon Y^M\to
Y$ the only map such that $f=\pi^Mf^M$. {\em If $f$~is fully
closed, then $\mathcal{M}$ is upper semicontinuous, $Y^M$~is a
regular space, and $f^M,\pi^M$ are fully closed maps.}

A proof of this lemma (cf.\ \cite[Lemma 1.2.9]{eng1}) is routine.

\begin{statem}\label{ext_partit}
{\bf Lemma.} Suppose that $M,A,B\subset X$ are closed subsets of a
normal space~$X$, $A\cap B=\emptyset$, and $L$ is a partition in
$M$ between $M\cap A$ and $M\cap B$. If $X\setminus L$ is a normal
space, then there are disjoint open sets $U,V\subset X$ such that
$A\subset U$, $B\subset V$, $M\setminus L=(U\cup V)\cap M$, and\/
$\cl U\cap\cl V\subset L$.\qed
\end{statem}

\begin{proof}[Proof of Theorem \ref{indor}]
We start with some general construction for arbitrary $X$, $f$,
$Y\!=\!fX$, and disjoint closed sets $A,B\!\subset\! X$. We can
assume that $p\!=\!\Indo Y\!<\!\infty$ and $q=\Indo f<\infty$.
Clearly\/ $p,q\geq 0$. Since $f$ is fully closed, $M=fA\cap fB$ is
a discrete closed subspace of $Y$. Consider $Y^M$, $f^M\colon X\to
Y^M$, and $\pi^M\colon Y^M\to Y$.  The restriction $f^M|f^{-1}M$
is a homeomorphism onto $N=(\pi^M)^{-1}M$, and we shall construct
a $G_\delta$ partition in~$Y^M$ between the disjoint sets $f^M A$
and $f^M B$. The pre-image $f^{-1}M$ is homeo\-mor\-ph\-ic to the
discrete sum of point-inverses $f^{-1}y$, $y\in M$, and hence,
$\Indo N\leq q$. There is a $G_\delta$ partition $L$ in $N$
between $N\cap f^M A$ and $N\cap f^M B$, where $\Indo L\leq q-1$.
As $f^M$ is a closed map, $Y^M$ is a normal space. Since $M\subset
Y$ is a $G_\delta$-set, $Y^M\setminus N$ and $Y^M\setminus L$ are
$F_\sigma$-sets in $Y^M$, and hence, they are also normal spaces
(see \cite[Exercise 2.1.E]{eng2}). By Lemma \ref{ext_partit},
there are disjoint open sets $U,V\subset Y^M$ such that $f^M
A\subset U$, $f^M B\subset V$, $N\setminus L=(U\cup V)\cap N$,
and\/ $\cl U\cap\cl V\subset L$. As $Y\setminus M$ is an open
$F_\sigma$ subset of~$Y$, it is a countable union of closed
$G_\delta$ subsets $F_i$ of $Y$. Since $\pi^M|Y^M\setminus N$ is a
homeo\-mor\-ph\-ism onto $Y\setminus M$, we obtain $\Indo
(Y^M\setminus N)\leq p$ by Theorem \ref{count}. Thus, there are
disjoint open sets $U',V'\subset Y^M\setminus N$ with $\cl
U\setminus N\subset U'$, $\cl V\setminus N\subset V'$, and
$L'=Y^M\setminus (N\cup U'\cup V')$ is a $G_\delta$-set with
$\Indo L'\leq p-1$. Observe that $L\cup L'=Y^M\setminus(U\cup
U'\cup V\cup V')$ is a $G_\delta$ partition in $Y^M$ between $f^M
A$ and $f^M B$. It follows that {\em $(f^M)^{-1}(L\cup L')$ is a
$G_\delta$ partition in $X$ between $A$ and~$B$.}

We proceed by induction on $p$. If $p=0$, then $L'=\emptyset$ and $\Indo (f^M)^{-1}(L\cup L')=\Indo L\leq q-1$. As we took
arbitrary sets $A$ and $B$, we have $\Indo X\leq q=p+q$. Assume the theorem is true for fully closed maps whose images have
$\Indo <p>0$. Then, $L'$ is the countable union of closed $G_\delta$-sets $L_i=L'\cap (\pi^M)^{-1}F_i\subset Y^M$ with
$\Indo L_i\leq p-1$. The restrictions $f^M|(f^M)^{-1}L_i\colon (f^M)^{-1}L_i\to L_i$ are fully closed, and by the induction
hypothesis we obtain $\Indo (f^M)^{-1} L_i\leq p+q-1$. Since $L$~and $(f^M)^{-1}L$ are homeo\-mor\-ph\-ic, we have $\Indo
(f^M)^{-1}(L\cup L')\leq p+q-1$ by Theorem \ref{count}. We have shown that $\Indo X\leq p+q$ because $(f^M)^{-1}(L\cup L')$
is a $G_\delta$ partition between disjoint closed sets $A$ and $B$, which were taken arbitrarily.
\end{proof}

\begin{statem}\label{ind}
{\bf Corollary.} Suppose that\/ $f$ is a fully closed map from a
non-empty normal space $X$ onto a perfectly normal space. If every
point-inverse of $f$ is perfectly normal, then\/ $\Ind X\leq\Ind
fX+\Ind f$.\qed
\end{statem}

The foregoing corollary may be considered as an $\Ind$-analogue of Fedorchuk's Theorem 4 in \cite{fed5}, which was stated
for $\ind$ and special, resolution fully closed maps~$f$. In a recent paper \cite[pp.\ 117--120]{fed_dg} Fedorchuk proves
the inequality for $\Ind$ and resolution maps $f$, where $fX$ are (metric, compact) two-manifolds.

Theorem \ref{indor} and Corollary \ref{ind} enable estimating
inductive dimensions in some well-known constructions (see
\cite{fed} for a survey). In particular, Fedorchuk's continua $B$
(\cite{fed5})---let us write $B_n$ instead---were the first
examples of {\em separable and first countable} compact spaces
with non-coinciding dimensions $\dim$ and $\ind$. Fedorchuk proved
that $\dim B_n=n$ and $2n-1\leq\ind B_n\leq 2n$. Since each $B_n$
has a fully closed map onto the $n$-dimensional sphere, and every
point-inverse of the map is homeomorphic to the $n$-dimensional
torus, we obtain

\begin{statem}\label{Indfed}
{\bf Corollary.} Fedorchuk's continua $B_n$ have also\/ $\Ind
B_n\leq\Indo B_n\leq 2n$.\qed
\end{statem}

\noindent (In fact, we shall see that $\Indo B_n=2n$ by Theorem
\ref{indok}.)

Chatyrko \cite{chat5} constructed separable first countable
continua $I_n$ and $(S_1)_n$, and proved that $I_n$ are chainable,
$(S_1)_n$ are ho\-mo\-ge\-ne\-ous, $\dim I_n=\dim (S_1)_n=1$, and
$\ind I_n=\ind (S_1)_n=n$.

\begin{statem}\label{Indchat}
{\bf Corollary.}  Chatyrko's continua $I_n$ and $(S_1)_n$ have
also
$$\Ind I_n=\Ind (S_1)_n=\Indo I_n=\Indo (S_1)_n=n.$$
\end{statem}

\begin{proof}
There is a sequence
$\mbox{}\ldots\stackrel{\pi^{n+1}_n}{\longrightarrow} I_n
\stackrel{\pi^n_{n-1}}{\longrightarrow}\ldots
\stackrel{\pi^3_2}{\longrightarrow}I_2
\stackrel{\pi^2_1}{\longrightarrow}I_1=[0,1]$ of fully closed onto
maps $\pi_n^{n+1}$, see \cite{chat5}. For each $n$ and every $t\in
I_n$, the pre-image $(\pi^{n+1}_n)^{-1}t$ is homeomorphic to
$[0,1]$. Using induction and Theorem \ref{indor}, we obtain $\Indo
I_n\leq n$.

In the case of $(S_1)_n$ there exists an analogous sequence of
maps, whose point-inverses are homeomorphic to a circumference.
\end{proof}

In Corollary \ref{ind} one can replace those perfectly normal
spaces by another class of spaces in which $\Ind=\Indo$. This
could be the class of hereditarily perfectly $\kappa$-normal
spaces (Fedorchuk \cite{fed2}); surely, every discrete closed
subset of $fX$ should be $G_\delta$ in $fX$.

The assumption that $f$ is fully closed is necessary in Corollary \ref{ind}. Under a set-theoretical assumption consistent
with $\mathsf{ZFC}$, Fedorchuk \cite{fed4} constructed a perfect map $f_{\rm F}\colon X_{\rm F}\to Y_{\rm F}$, where $X_{\rm
F}$ and $Y_{\rm F}$ are perfectly normal, locally compact, and countably compact spaces, $\dim X_{\rm F}\!=\!\Ind X_{\rm
F}\!=\!1$, $\Ind Y_{\rm F}\!=\!0$, and $\Ind f_{\rm F}^{-1}y\!=\!0$ for every $y\in Y_{\rm F}$. On the other hand, it is not
sufficient to assume only that $f$ is fully closed. Chatyrko \cite{chat6} has constructed a certain fully closed map $f_{\rm
Ch} \colon\! X_{\rm Ch}\!\to\! A_\mathfrak{c}$ from a compact space $X_{\rm Ch}$ with $\Ind X_{\rm Ch}\!=\!\Indo X_{\rm
Ch}\!=\!2$ onto the compact space~$A_\mathfrak{c}$ with the only accumulation point~$y_0$, $\card
A_\mathfrak{c}=\mathfrak{c}$. All point-inverses $f_{\rm Ch}^{-1}y$, where $y_0\neq y\in A_\mathfrak{c}$, are single points,
and $1=\Ind f_{\rm Ch}^{-1}y_0<\Indo f_{\rm Ch}^{-1}y_0=2$. For some other maps $f\colon X\to Y$, even $\Ind X-\Ind fX-\Ind
f>1$. For every pair of natural numbers $m>n\geq 1$, the present author \cite{krz} has constructed a compact space $X_{m,n}$
such that $\Ind X_{m,n}=m$ and every component of $X_{m,n}$ is homeo\-morphic to the $n$-dimensional cube $[0,1]^n$. In
consequence, if $\mathcal{D}$~stands for the decomposition of $X_{m,n}$ into its components, and $f_{\rm K}\colon X_{m,n}\to
X_{m,n}/\mathcal{D}$ is the natural quotient\linebreak map (it is not fully closed), then $\Ind X_{m,n}-\Ind
X_{m,n}/\mathcal{D}-\Ind f_{\rm K}=m-n$.

\begin{statem}\label{full}
{\bf Lemma.} Suppose that $f$ is a fully closed map from a normal
space $X$, $L\subset X$ is a closed $G_\delta$-set, and
$A,B\subset fX$ are disjoint closed sets. Then,
\begin{wylicz}
\item $fL\cap
f(X\setminus L)$ is the countable union of discrete closed
subspaces of $fX$ {\rm (cf.\ \cite[Definition 3 and Lemma
2]{fed1})}.
\end{wylicz}
If moreover every discrete closed subspace of $fX$ is a
$G_\delta$-subset, then
\begin{wylicz}\setcounter{licz}{1}
\item
$fL$ is a $G_\delta$-set in $fX$; and
\item whenever $L$ is a partition between $f^{-1}A$ and $f^{-1}B$,
there is a countable family of discrete closed sets
$\Gamma_i\subset fX\setminus(A\cup B)$ such that the union
$fL\cup\bigcup_{\,i} \Gamma_i$ is a $G_\delta$ partition in $fX$
between $A$ and $B$.
\end{wylicz}
\end{statem}

\begin{proof}
(a) There is a sequence of closed sets $F_i\subset X$ such that
$X\setminus L=\bigcup_{\,i} F_i$. The intersections $fF_i\cap fL$
are discrete and closed, and $f(X\setminus L)\cap
fL=\bigcup_{\,i}(fF_i\cap fL)$. (b) If $fF_i\cap fL$ are
$G_\delta$ in $fX$, then $fF_i\setminus fL$ are $F_\sigma$. Since
$fX\setminus fL=\bigcup_{\,i}(fF_i\setminus fL)$, $fL$ is
$G_\delta$. The assertion (c) appears to be implicitly shown in
the proof of \cite[Theorem III.2.6]{fed} if one can use (b). The
new point is to apply Lemma \ref{full}(b) and prove that, if $P=L$
in \cite[p.\ 4247]{fed} is $G_\delta$, then $fP$, $U_i\cup
f^{-1}fP$, $f(U_i\cup f^{-1}fP)$, and $K=f(U_1\cup f^{-1}fP)\cap
f(U_2\cup f^{-1}fP)=fP\cup\bigcup_{j,k} \Gamma_{jk}$ are
$G_\delta$-sets.
\end{proof}

Applying induction, Lemma \ref{full}(c), Theorem \ref{count}, and
Proposition \ref{at}, we obtain the following two theorems. (The
first one is an $\Indo$-analogue of Theorem III.2.6 on $\Ind$ in
\cite{fed}.)

\begin{statem}\label{rlfc1}
{\bf Theorem.} If $f$ is a fully closed map from a normal space
$X$ to a space whose every discrete closed subspace is a
$G_\delta$-set, then\/ $\Indo fX\leq \Indo X+1$.\qed
\end{statem}

\begin{statem}\label{rlfc}
{\bf Theorem.} If $f$ is a ring-like fully closed map from a
compact space $X$ to a first countable space, then\/ $\Indo
fX\leq\Indo X$.\qed
\end{statem}

\begin{statem}\label{gie}
{\bf Lemma.} If $f$ is a ring-like fully closed map from a compact
space~$X$ onto a non-degenerate continuum, then every $G_\delta$
partition in $X$ contains a point-inverse of $f$.
\end{statem}

\begin{proof}
Take a $G_\delta$ partition $L$ in $X$. There are two cases. (1) If $fL$ is uncountable, then by Lemma \ref{full}(a),
$fL\cap f(X\setminus L)$ is countable, $fL\setminus f(X\setminus L)\ni y$ for some $y\in fX$, and $f^{-1}y\subset L$. (2)
Suppose that $fL$ is countable. $L$ contains a {\em thin} partition $F$, i.e. such that there are disjoint non-empty open
sets $U_1,U_2\subset X$ with $X\setminus F=U_1 \cup U_2$ and $F=\fr U_1=\fr U_2$. By the remark \ref{irr} and Proposition
\ref{at}, $f$ is a monotone irreducible map, and $X$~is a continuum. By the same argument as in Fedorchuk \cite[the proof of
Lemma 4, p.\ 167]{fed1}, we infer that $f^{-1}x\subset F$ for every point $x\in fX$ isolated in $fF$.
\end{proof}

\begin{statem}\label{indok}
{\bf Theorem.} If $f$ is a ring-like fully closed map from a
non-empty compact space $X$ onto a first countable space $Y$,
and\/ $\Indo f^{-1}y=\Indo f$ for every $y\in Y$, then $\Indo X=
\Indo Y + \Indo f$.
\end{statem}

\begin{proof}
Theorem \ref{indor} yields the inequality ,,$\leq$''. Fix $m=\Indo f$. The inequality ,,$\geq$'' will be proved by induction
on $n=\Indo X<\infty$. Clearly, $n\geq m$. Let $n=m$, and suppose that $\Indo Y>0$. Then, $Y$~has a non-degenerate
component~$Y'$. By Lemma \ref{gie}, every $G_\delta$ partition $L$ in $X'=f^{-1}Y'$ contains a point-inverse of~$f$ and has
$\Indo L\geq m$. Hence, $\Indo X'>n$. A contradiction. Thus, $\Indo Y=0$ and the inequality ,,$\geq$'' is true.

Assume $n>m$. In order to show that $\Indo Y\leq n-m$, take a pair
of disjoint closed sets $A,B\subset Y$. There is a $G_\delta$
partition $L$ in $X$ between $f^{-1}A$ and $f^{-1}B$, $\Indo L\leq
n-1$. It is a consequence of Lemma \ref{full}(c) that there is a
$G_\delta$ partition $K\supset fL$ in $Y$ between $A$ and~$B$,
where $K\setminus fL$ is countable. Lemma \ref{full}(a) implies
that also $K\cap f(X\setminus L)=\{y_i\colon i=1,2,\ldots\}$. We
have $f^{-1}K=L\cup\bigcup_{\,i}f^{-1}y_i$, and $\Indo f^{-1}K\leq
n-1$ by Theorem \ref{count}. By the induction hypothesis, $\Indo
K\leq n-m-1$. We have shown that $\Indo Y\leq n-m$.
\end{proof}

The following statement may be considered as a generalization of Theorem 3 in Fedorchuk \cite{fed5}.

\begin{statem}\label{cm}
{\bf Theorem.} Suppose that $f\colon X\to Y$ is a surjective ring-like map between compact spaces~$X$ and $Y$, and\/ $\dim
Y\geq 1$. If\/ $\ind f^{-1}y\geq m$ for every $y\in Y$, then\/ $\ind X\geq\dim Y+m-1$.
\end{statem}

\begin{proof}
If $\dim Y=\infty$, then $\ind X\geq\dim X=\infty$ by Proposition \ref{ch1} (and \cite[Theorem 3.1.29]{eng1}). We can assume
that $\dim Y<\infty$. It suffices to prove by induction that for every natural number $k\geq1$, the inequality\/ $\dim Y\geq
k$ implies $\ind X\geq k+m-1$. For $k=1$ the implication is obvious.

We shall use the following classical notion. A compact space $M$
with $\dim M\!=\!n$ is called an {\em $n$-dimensional Cantor
manifold}\/ provided that every partition $L$ in~$M$ has $\dim
L\geq n-1$. P.S.~Alexandroff \cite{al} proved that every compact
space $Z$ with $1\leq n=\dim Z<\infty$ contains an $n$-dimensional
Cantor manifold (see also \cite[Theorem 1.9.9 and Exercise
3.2.F]{eng1}).

Let $n=\dim Y\geq k\geq 2$. We can assume that $Y$ is an $n$-dimen\-sion\-al Cantor manifold. Then, $X$~is a continuum and
$f$ is an irreducible monotone map by \ref{irr} and Proposition \ref{at}. If $L$ is a partition in $X$, then $fL$ is a
partition in~$Y$ (Proposition \ref{folk}), and $\dim fL\geq n-1\geq k-1$. Now, $fL$ contains a component $P$ with $\dim
P\geq n-1$. Since $P=f(L\cap f^{-1}P)$, the remark \ref{irr} implies $f^{-1}P\subset L$. By the obvious induction
hypothesis, $\ind L\geq\ind f^{-1}P\geq k+m-2$. Therefore, $\ind X\geq k+m-1$.
\end{proof}

\section{Main tools for constructions}

For our constructions we need ingredients of two types. The
following resolution theorem (the first type) is a modification of
well-known results.

\begin{statem}\label{main}
{\bf \hspace{-.1ex}Theorem } \em \hspace{-.1ex}(cf.\
\hspace{-.1ex}Chatyrko \hspace{-.1ex}\cite{chat1},
\hspace{-.1ex}Emeryk \hspace{-.1ex}\cite{em1},
\hspace{-.1ex}Fedorchuk \hspace{-.1ex}\cite{fed5}). \em
\hspace{-.1ex}Suppose \hspace{-.1ex}that $X$~\hspace{-.1ex}is a
first countable continuum, and for every\/ $x\in X$, $Y_x$ is a
metrizable continuum. Then, there exists a first countable
continuum\/ $Z\!=\!Z(X,Y_x)$ with a map $\pi\colon Z\to X$ such
that
\begin{wylicz}
\item for every $x\in X$, the pre-image
$\pi^{-1}x$ is homeomorphic to $Y_x$; and
\item $\pi$ is ring-like and fully closed.
\end{wylicz}

Moreover, the conjunction of\/ {\em (a)} and\/ {\em (b)} implies
that
\begin{wylicz}\setcounter{licz}{2}
\item $\dim Z=\max\{\dim X,\dim\pi\}$;
\item if $X$ is perfectly normal, then\/ $\Ind Z\leq\Ind X+\Ind\pi$;
\item if $X$ is separable, then so is $Z$;
\item if all $Y_x$ are non-degenerate continua, and $P\subset Z$ is a
metrizable continuum, then the image $\pi P$ is a single point;
and
\item if $X$ and all\/ $Y_x$ are hereditarily indecomposable\/
{\em [}hereditarily decomposable\hspace{.2ex}{\em ]}, then so is
$Z$.
\end{wylicz}
\end{statem}

\begin{proof}[Sketch of proof]
Since Chatyrko's paper \cite{chat1} has not been translated into
English, we sketch his construction for the convenience of the
reader\footnote{In a forthcoming paper, joint with M.G.\
Charalambous, we describe a generalization of Chatyrko's
construction in more detail.}. We can assume that each $Y_x$ is a
subspace of the Hilbert cube $[0,1]^\infty$. Using the local
connectedness of $[0,1]^\infty$, one constructs a map
$g_x\colon(0,1]\to [0,1]^\infty$ such that for each natural number
$n$, $Y_x\subset\cl g_x (0,1/n]\subset\ball{Y_x}{1/n}$, where B
stands for a ball. One takes a map $f_x\colon X\to [0,1]$ with
$f^{-1}0=\{x\}$, and writes $h_x\colon X\setminus\{x\}\to
[0,1]^\infty$ for the composition $g_x(f_x|X\setminus\{x\})$. $Z$
is the set $\bigcup\{\{x\}\times Y_x\colon x\in X\}\subset X\times
[0,1]^\infty$, and $\pi\colon Z\to X$, $\pi (x,y)=x$. The topology
on $Z$ is generated by a base of neighborhoods at any point
$(x,y)\in Z$; the base consists of all sets
$$W(V,U)=[\{x\}\times (V\cap Y_x)]\cup \pi^{-1}(U\cap h_x^{-1}V),$$
where $U\subset X$ and $V\subset [0,1]^\infty$ are neighborhoods of $x\in X$ and $y\in Y_x$ respectively. The above is a
generalization of Fedorchuk's construction \cite{fed5} (cf.\ also Fedorchuk \cite[the proof of Lemma 1]{fed1} and
\cite[Section III.1]{fed}). One checks that $Z$ is a first countable continuum, and $\pi$ satisfies (a) and (b).

Corollary \ref{rdim} yields the equality (c), and Corollary \ref{ind} yields (d). The statement (e) is an easy consequence
of the fact that $\pi$ is an irreducible map (see the remark \ref{irr}). The statement (f) follows from \ref{irr} and
Proposition \ref{met}, and (g)~is a simple property of atomic maps.
\end{proof}

\begin{statem}\label{r1}
{\bf Remarks.} \em (1) If a subcontinuum $P$ of $Z=Z(X,Y_x)$ in
Theorem \ref{main} is non-metrizable, then $P=\pi^{-1}\pi P$ by
\ref{irr}. Hence, {\em if all $Y_x$ are non-degenerate continua,
then for every non-degenerate continuum $P\subset Z$ and any point
$z\in P$, there is a non-degenerate metrizable continuum $Q$ such
that $z\in Q\subset P$.}

(2) One can combine the proofs by Chatyrko \cite{chat1} and
Fedorchuk \cite[Lemma~1]{fed1} in order to obtain a map $\pi$ that
satisfies also the assertion (2) of Lemma 1 in \cite{fed1}. This
enables one to construct the continuum $Z=Z(X,Y_x)$ under the
continuum hypothesis (see \cite[pp.\ 166-167]{fed1}) so that
\begin{wylicz}\item[(\dag)\hfill]
\em if\/ $\mathsf{CH}$ is true, and the continuum $X$ in Theorem
\ref{main} is perfectly normal and hereditarily separable, then
$Z$ is perfectly normal and hereditarily separable.
\end{wylicz}
\end{statem}

We shall also need Cook continua (the second type of ingredients),
whose sub\-continua will be taken as the $Y_x$'s of
Theorem~\ref{main}.

\begin{statem}\label{m1}\em
{\bf Example} (Cook \cite{cook}; see also A.\ Pultr, V.\
Trnkov\'{a} \cite[Appendix A]{pul} for a detailed construction).
\em There exists a metrizable, one-dimensional, hereditarily
indecomposable Cook continuum $M_1$ that does not contain
non-degenerate continuous images of plane continua.
\end{statem}

\begin{proof}
J.W.\ Rogers, Jr.\ \cite{rog} observed that Cook's continuum
\cite{cook} does not contain non-degenerate continuous images of
plane continua.
\end{proof}

\begin{statem}\label{ma}
{\bf Example} \em (Maćkowiak \cite{mac15}). \em There exists a
metrizable, chainable, hereditarily decomposable Cook
continuum.\qed
\end{statem}

\section{Anderson-Choquet (and similar) continua with $\dim >1$}
\label{ach}

The following continua are neither Anderson-Choquet nor Cook, but
they are HI analogues of Fedorchuk's spaces with non-coinciding
dimensions (\cite{fed5}, cf.\ also our Corollary \ref{Indfed}).

\begin{statem}
{\bf Theorem.} For every natural number $n\geq 1$, there exists a
non-metrizable, separable, first countable, hereditarily
indecomposable continuum $Z$ such that\/ $\dim Z=n$~and
$2n-1\leq\ind Z\leq\Ind Z\leq\Indo Z=2n$.
\end{statem}

\begin{proof}
By a theorem of R.H.\ Bing \cite{bing}, there exists a metric HI continuum $X$ with $\dim X=n$. Let $Y_x=X$ for $x\in X$,
apply Theorem \ref{main}, and put $Z=Z(X,Y_x)$.

The properties of $Z$ are consequences of the statements \ref{indok}, \ref{cm}, and \ref{main}(c--g).
\end{proof}

Let us notice that the first examples of non-metrizable HI
continua were constructed by Emeryk \cite{em2}.

\begin{statem}\label{ach_dim}
{\bf Theorem.} For every natural number $n\geq 1$, there exists a
non-metrizable, separable, first countable, hereditarily
indecomposable Anderson-Choquet continuum $Z$ with $n=\dim
Z\leq\ind Z\leq\Ind Z\leq\Indo Z=n+1$.
\end{statem}

\begin{proof}
Let $X$ be a metric $n$-dimensional HI continuum (Bing
\cite{bing}). Take a family $\{C_x\colon x\in X\}$ of pairwise
disjoint non-degenerate subcontinua of Cook's continuum $M_1$
(Example \ref{m1}), apply Theorem \ref{main}, and put
$Z=Z(X,C_x)$.

Most of the desired properties of $Z$ follow from Theorems \ref{indok} and \ref{main}(c-g). There remains to prove that $Z$
is an Anderson-Choquet continuum. Let $\pi\colon\! Z\!\to\! X$ be the map of Theorem \ref{main}, and choose any
non-degenerate continuum $P\subset Z$ and any embedding $\varphi\colon P\to Z$. Take an arbitrary point $z\in P$. By Remark
\ref{r1}(1), $P$~contains a non-degenerate metrizable continuum $Q\ni z$. The statement (f) of Theorem \ref{main} guarantees
that $Q\subset\pi^{-1}x$ and $\varphi Q\subset\pi^{-1}x'$ for some $x,x'\in X$. Since $\pi^{-1}x$ and $\pi^{-1}x'$ are
homeo\-morph\-ic to $C_x$ and $C_{x'}$, respectively, we obtain $x=x'$, $\varphi|Q=\id_Q$, and $\varphi z=z$. We have shown
that $\varphi=\id_P$.
\end{proof}

We shall adapt the foregoing construction and proof in order to
obtain hereditarily rigid finite-group actions. We start with some
terminology. Let $X$ be a space, and $G$ a finite group. We write
$H(X)$ for the group of all homeomorphisms $X\to X$. Every
homomorphism $\xi\colon G\to H(X)$ is called a\/ {\em $G$-action
on $X$}; the value of $\xi$ at $g\in G$ will be denoted by
$g^\xi\in H(X)$. This $G$-action is said to be {\em
fixed-point-free} if for each $g\in G\setminus\{e\}$, the
homeo\-morph\-ism $g^\xi\colon X\to X$ does not have a fixed
point. Let $\zeta\colon G\to H(Y)$ be a $G$-action on a space $Y$.
A~map $f\colon X\to Y$ is said to be {\em equivariant}\/ if
$g^\zeta f=fg^\xi$ for each $g\in G$.

\begin{statem}\label{act}
{\bf Theorem.} Suppose that $X$ is a first countable continuum,
$G$ a finite group, and $\xi$ is a fixed-point-free $G$-action on
$X$. Then, there exists a first countable continuum $Z$ with a
fixed-point-free isomorphic $G$-action\/ $\zeta\colon G\to H(Z)$
and an equi\-var\-i\-ant map $\pi\colon Z\to X$ such that
\begin{wylicz}
\item $\pi$ is a ring-like fully closed map, and all point-inverses of
$\pi$ are metrizable one-dimensional continua;
\item for every non-degenerate continuum $P\subset Z$ and every
embedding $\varphi\colon P\to Z$, there is a homeomorphism
$g^\zeta\in H(Z)$ such that $g^\zeta|P=\varphi$.
\end{wylicz}
\end{statem}

An important point in this theorem is that $\dim Z=\dim X$ by
Corollary \ref{rdim}, and $\Indo Z=\Indo X+1$ by Theorem
\ref{indok}. In the proof it will be seen that we can ensure that
all point-inverses of $\pi$ are HI (if we use subcontinua of
Example \ref{m1} in our construction), or alternatively, that they
are HD (if we use Example~\ref{ma}).

We may apply this theorem to some standard examples of group
actions. It is well-known that for every finite group $G$ and
every $n\geq 2$, there exists a (compact metric) $n$-manifold
without boundary with a fixed-point-free $G$-action (see J.\ de
Groot, R.J.\ Wille \cite[p.\ 444]{gw}). In case $n=1$, there
exists a connected finite graph with a fixed-point-free $G$-action
($G$ acts on its Cayley graph). Two more simple examples: Using
Anderson and Choquet's original HD continuum \cite{and}, a Cayley
graph of the group $G$, and the method from \cite{gw}, one easily
constructs a metric one-dimensional HD continuum $Z$ with a
fixed-point-free $G$-action $\zeta$ that is an isomorphism $G\to
H(Z)$ and satisfies the assertion (b) of Theorem~\ref{act}. Using
our Anderson-Choquet continuum (of Theorem \ref{ach_dim}) instead
(and the same Cayley graph method), one constructs a
non-metrizable, separable, and first countable continuum $Z$ with
$\dim Z=n$ and a similar $G$-action $\zeta$ on $Z$.

Let us notice that there are numerous papers on group
representations in topology (see de Groot \cite{gr} for example)
and in the more general context of category theory (see the
bibliography in Pultr and Trnkov\'{a} \cite{pul}).

\begin{proof}[Proof of Theorem \ref{act}.] (I) Consider the family
$\mathcal{D}$ of all orbits $\{g^\xi x\colon\! g\!\in\! G\}$,
$x\!\in\! X$, and the quotient space $X'=X/\mathcal{D}$. The
quotient projection $q_\xi\colon X\to X'$ is a covering map. Take
a family $\{C_{x'}\colon x'\in X'\}$ of pairwise disjoint
non-degenerate subcontinua of Cook's continuum $M_1$ (Example
\ref{m1}). Use Theorem \ref{main}, and take $Z'=Z(X',C_{x'})$ and
$\pi'\colon Z'\to X'$, a map that satisfies the statements (a--g)
of Theorem \ref{main}. Consider the set $Z=\bigcup_{x'\in X'}
(q_\xi^{-1}x'\times\pi'^{-1} x')\subset X\times Z'$. We define
$\pi(x,t)=x$, $q_\zeta (x,t)=t$, and $g^\zeta(x,t)=(g^\xi x,t)$
for $x\in X$, $t\in\pi'^{-1}q_\xi x$, and $g\in G$. The following
diagram commutes,
\[\mbox{\begin{picture}(180,95)
\put(75,87){$X$}\put(159,87){$X'$} \put(75,34){$Z$}\put(159,34){$Z'$}
\put(23,53){$X$}\put(23,0){$Z$}\put(36,57){\vector(4,1){119}} \put(36,4){\vector(4,1){119}}\put(72,33){\vector(-3,-2){36}}
\put(72,86){\vector(-3,-2){36}}\put(87,38){\vector(1,0){68}} \put(87,91){\vector(1,0){68}}\put(80,47){\vector(0,1){34}}
\put(164,47){\vector(0,1){34}}\put(28,13){\vector(0,1){34}}

\put(115,95){$\scriptstyle q_\xi$}\put(115,42){$\scriptstyle
q_\zeta$}\put(20,22){$\scriptstyle \pi$}\put(73,51){$\scriptstyle
\pi$}\put(154,57){$\scriptstyle \pi'$}\put(112,69){$\scriptstyle
q_\xi$}\put(112,16){$\scriptstyle
q_\zeta$}\put(49,79){$\scriptstyle
g^\xi$}\put(49,26){$\scriptstyle g^\zeta$}
\end{picture}}\]
and $\zeta$, $g\mapsto g^\zeta$, is a monomorphism from $G$ to the
group of permutations of $Z$. Finally, we equip $Z$ with the
smallest topology such that $\pi$ and $q_\zeta$ are continuous.
Observe that, if $U\subset X$ is an open set such that
$q_\xi|U\colon U\to q_\xi U$ is a homeomorphism, then
$q_\zeta|\pi^{-1} U\colon \pi^{-1} U\to \pi'^{-1}q_\xi U$ is also
a homeomorphism. Indeed, $q_\zeta|\pi^{-1} U$ is one-to-one, all
open subsets of $\pi^{-1}U$ have the form
$q_\zeta^{-1}V\cap\pi^{-1}U$, where $V\subset Z'$ are open, and
$q_\zeta(q_\zeta^{-1}V\cap\pi^{-1}U)=V\cap
q_\zeta\pi^{-1}U=V\cap\pi'^{-1}q_\xi U$. It follows that $q_\zeta$
is a closed covering map, and hence, $Z$ is a compact space by
\cite[Theorem 3.7.2]{eng2}. Similarly, $g^\zeta$ are
homeomorphisms of $Z$. In view of the above diagram, $\pi$ is an
equivariant map.

(II) Every point-inverse $\pi^{-1}x$ is homeomorphic to
$q_\zeta\pi^{-1}x=\pi'^{-1}q_\xi x$ and $C_{q_\xi x}$. As $\pi$ is
a monotone map, the connectedness of $X$ implies the connectedness
of~$Z$. Using the fact that $\pi'$ is ring-like, one easily checks
that also $\pi$ is ring-like.

In order to prove that $\pi$ is fully closed take disjoint closed
sets $A,B\subset Z$. By virtue of the compactness of $X$, it is
sufficient to show that, {\em if\/  $\cl U\subset X$ and\/
$q_\xi|\cl U$ is one-to-one, then\/ $\cl U$ contains a finite
number of points in $\pi A\cap\pi B$.} Indeed, consider the set
$F=\pi^{-1}\cl U$ and the restriction $\pi|F$. The sets $g^\zeta
F$, $g\in G$, are pairwise disjoint, and $q_\zeta|F$ is a
one-to-one function. Thus, $q_\zeta (F\cap A)$ and $q_\zeta (F\cap
B)$ are disjoint closed subsets of $Z'$, and the intersection
$$\pi'q_\zeta (F\cap
A)\cap\pi' q_\zeta(F\cap B)=q_\xi\pi(F\cap A)\cap q_\xi\pi(F\cap B)= q_\xi(\cl U\cap \pi A\cap\pi B)$$ is finite since
$\pi'$ is a fully closed map. Therefore, $\cl U\cap \pi A\cap\pi B$ is finite, and $\pi$ is fully closed. Moreover, $\pi$
satisfies statements analogous to Theorem \ref{main}(c--g) and Remark \ref{r1}(1).

(III) Let $P\subset Z$ be a non-degenerate continuum, and
$\varphi\colon P\to Z$ an embedding. We claim that {\em for every
$z\in P$ there is a $g\in G$ such that $\varphi z=g^\zeta z$.}
Indeed, by Remark \ref{r1}(1), $P$ contains a non-degenerate
metrizable continuum $Q\ni z$. By a statement analogous to Theorem
\ref{main}(f), there are $x,x'\in X$ such that $Q\subset\pi^{-1}x$
and $\varphi Q\subset\pi^{-1}x'$. The restrictions $q_\zeta|Q$ and
$q_\zeta\varphi|Q$ are embeddings into $\pi'^{-1}q_\xi x$ and
$\pi'^{-1}q_\xi x'$, respectively. Since these point-inverses are
homeomorphic to $C_{q_\xi x}$ and $C_{q_\xi x'}$, respectively, we
have $q_\xi x=q_\xi x'$, $q_\zeta|Q=q_\zeta\varphi|Q$, and
$q_\zeta z=q_\zeta\varphi z$. Hence, there is a $g\in G$ such that
$\varphi z=g^\zeta z$. The foregoing claim implies that
$P=\bigcup_{g\in G} F_g$, where $F_g=\{z\in P\colon \varphi
z=g^\zeta z\}$. $F_g$ are closed, pairwise disjoint, and $P$ is
connected. In consequence, only one $F_g$ is non-empty. Thus,
there is a $g\in G$ such that $P=F_g$ and $\varphi=g^\zeta|P$.

When we apply the statement (b) to $P=Z$, we infer that $\zeta$ is
an isomorphism onto $H(Z)$.
\end{proof}

A non-degenerate continuum~$X$ will be called a {\em weak Cook
continuum\footnote{There is some difference in~terminology: in
[\ref{ma1}--\ref{ma2}] our weak Cook continua are just called Cook
continua.}}\/ if for every subcontinuum $P$ of $X$, every map
$f\colon P\to X$ with $P\cap fP=\emptyset$\/ is constant.

\begin{statem}\label{retr}
{\bf Proposition} \em (Maćkowiak \cite[Proposition
29]{mac1}\footnote{The proof of Proposition 29(i) in \cite{mac1}
works for arbitrary Hausdorff continua.}). \em Suppose that $X$ is
a weak Cook continuum. If $P$ is a subcontinuum of $X$ and
$f\colon P\to X$ is a non-constant map, then $fP\subset P$ and $f$
is a monotone retraction.\qed
\end{statem}

\begin{statem}\label{kuka}
{\bf Theorem.} There exists a non-metrizable, separable, first
countable Cook continuum $Z$ with $2=\dim Z\leq\ind Z\leq\Ind
Z\leq\Indo Z=3$.
\end{statem}

\begin{proof}
Take the square $[0,1]^2$ and a family $\{C_x\colon x\in
[0,1]^2\}$ of pairwise disjoint non-degenerate subcontinua of
Cook's continuum $M_1$. Put $Z=Z([0,1]^2,C_x)$, and let $\pi\colon
Z\to [0,1]^2$ be the map of Theorem \ref{main}.

Most of the desired properties of $Z$ follow from Theorems \ref{indok} and \ref{main}(c-g). We shall prove that $Z$ is a
weak Cook continuum. Assume that $P$ is a subcontinuum of~$Z$, and $f\colon P\to Z$ is a map with $P\cap fP=\emptyset$. We
claim that {\em for every metrizable continuum $Q\subset P$, the restriction $f|Q$ is constant.} Indeed, $fQ$ is a
metriz\-able subcontinuum of $Z$. Theorem \ref{main}(f) implies that $Q$ and $fQ$ are homeomorphic to disjoint subcontinua
of Cook's $M_1$. Hence, $fQ$ is a single~point. Thus, if $P$ is metrizable, we are done. If not, $P=\pi^{-1}\pi P$ by
\ref{irr}, and the above claim implies that there is a factorization $f=g(\pi|P)$, where $g\colon\pi P\to Z$ is continuous.
Hence, $fP$ is a metrizable continuum contained in a point-inverse of $\pi$. Since $M_1$ contains only degenerate images of
plane continua, $g$ and $f$ are constant maps.

Now, choose a continuum $P\!\subset\! Z$ and a non-constant map $f\colon \!P\!\to\! Z$. Suppose {\em a contrario} that
$f\neq\id_P$. Take a point $z\in P$ with $fz\neq z$. It is a consequence of Proposition \ref{retr} that $A=f^{-1}fz$ and
$B=fP$ are non-degenerate continua with $A\cap B=\{fz\}$. Since $A$ is a retract of $A\cup B\ni z$ and $\pi^{-1}\pi z$ is a
Cook continuum, $\pi (A\cup B)$ is not a single point. Hence, $A\cup B=\pi^{-1}\pi (A\cup B)=\pi^{-1}\pi A\cup \pi^{-1}\pi
B$ by \ref{irr}, $\pi A\not\subset\pi B$, $\pi B\not\subset \pi A$, $A=\pi^{-1}\pi A$, and $B=\pi^{-1}\pi B$. Thus,
$\pi^{-1}\pi fz\subset A\cap B$. A contradiction. Therefore, $f=\id_P$, and $Z$~is a Cook continuum.
\end{proof}

\section{Chainable Cook continua with $\ind >1$}
\label{secchain}

Let $X$ be a chainable continuum. Elements $a\neq b$ of $X$ are
called {\em opposite end points}\/ if every open cover of $X$ has
a finite closed refinement $F_1,F_2,\ldots,F_k$~such that $a\in
F_1$, $b\in F_k$, and $F_i\cap F_j\neq\emptyset$ iff $|i-j|\leq
1$. Bing \cite[Theorem 15]{bing2} proved that {\em every
non-degenerate metric chainable continuum contains a chainable
continuum with a pair of opposite end points.}

The third ingredient for our construction in this section is the
following series of chainable continua, which were mentioned in
Corollary \ref{Indchat} (now, we need more detail).

\begin{statem}\label{snake}
{\bf Example} \em (Chatyrko \cite{chat5}; see also \cite{chat2,cf}
for $n\!=\!2,3$). \em There exists an inverse sequence
$(I_n,\pi_m^n)_{n,m=1}^\infty$, where $I_n$ are separable, first
countable, hereditarily decomposable chainable continua, and
$\pi_m^n\colon I_n\to I_m$ are surjective maps such that
\begin{wylicz}
\item $I_1$ is homeomorphic to\/ $[0,1]$;
\item each $\pi^{n+1}_n$ is atomic and fully closed, and each $\pi^n_1$ is
ring-like;
\item for each $n>1$, $m\in\{1,n-1\}$, and every $t\in I_m$, the
pre-image $(\pi^n_m)^{-1}t$ is homeomorphic to $I_{n-m}$;
\item for each $n$ and every pair\/ $0\leq s< t\leq 1$, the
pre-image $(\pi^n_1)^{-1}[s,t]$ is a chainable continuum with
opposite end points $a_{s,t}\in (\pi^n_1)^{-1}s, b_{s,t}\in
(\pi^n_1)^{-1}t$;
\item  for each $n$, every partition in~$I_n$ contains some point-inverse
$(\pi^n_1)^{-1}t$, where $t\in I_1$; and
\item $\ind I_n=\Ind I_n=\Indo I_n=n$ for each $n$ {\em (see our Corollary \ref{Indchat})}.\qed
\end{wylicz}
\end{statem}

\begin{statem}\label{chain}
{\bf Lemma} \em (Chatyrko \cite[Lemma 1]{chat5}). \em Suppose that
$X, Y$ are compact spaces, and $f\colon X\to Y$, $g\colon Y\to
[0,1]$ are surjective maps such that
\begin{wylicz}
\item $f$ is fully closed, and both $g$ and $gf$ are ring-like;
\item for every\/ $t\in [0,1]$, the pre-image $(gf)^{-1}t$ is a
chainable continuum with a pair of opposite end points; and
\item for every pair\/ $0\leq s< t\leq 1$, the pre-image
$g^{-1}[s,t]$ is a chainable continuum with opposite end points
$a_{s,t}\in g^{-1}s, b_{s,t}\in g^{-1}t$.
\end{wylicz}
Then, $X$ is a chainable continuum with a pair of opposite end
points $a\in (gf)^{-1}0$, $b\in (gf)^{-1}1$.\qed
\end{statem}

\begin{statem}\label{cooksnake}
{\bf Theorem.} For every natural number\/ $n\geq 1$, there exists a non-metrizable, separable, first countable, chainable,
hereditarily decomposable Cook continuum~$Z$ such that $n\leq\ind Z\leq\Ind Z\leq\Indo Z=n+1$, and every partition~$L$ in
$Z$ has\/ $\ind L\geq n-1$.
\end{statem}

\begin{proof}
Let $M$ be Maćkowiak's Cook continuum of Example \ref{ma}, and
$I_n$ Chatyrko's continuum of Example \ref{snake}. We claim that
{\em $M$ contains an uncountable family of pairwise disjoint
non-degenerate subcontinua $M_x$, $x\in I_n$.} Indeed, $M$~is HD,
and by a theorem of Bing \cite[Theorem 8]{bing2}, there is a
monotone surjective map $f\colon\! M\!\to\! [0,1]$. If
$0\!<\!t\!<\!1$ and the point-inverse $f^{-1}t$ were a single
point, $f^{-1}[0,t]$ would be a retract of $M$, and $M$ would not
be a Cook con\-ti\-n\-uum. Hence, $f^{-1}t$ is a non-degenerate
continuum if $0<t<1$. As $\card I_n=2^{\aleph_0}$, the claim is
proved\footnote{Recall that a continuum $X$ is said to be {\em
Suslinian}\/ if every family of pairwise disjoint non-degenerate
subcontinua of $X$ is countable. By a similar argument, we infer
that {\em no metric Cook continuum $X$ is Suslinian.} If $X$
contains an indecomposable continuum, then it is obviously not
Suslinian (see for instance \cite[Lemma 5.5]{krz1}). If $X$ is HD,
then we can assume that it is irreducible by \cite[\S 48 I,
Theorem 1]{kur}, and again there is a surjective monotone map
$f\colon X\to [0,1]$ by Kuratowski's theorems \cite[pp.\ 200 and
216]{kur}.

On the other hand, Maćkowiak \cite[Theorem 30]{mac1} constructed
an example of a metrizable, Suslinian, chainable, weak Cook
continuum.}. By \cite[Theorem 15]{bing2}, we can assume that every
$M_x$ has a pair of opposite end points. Finally, we apply Theorem
\ref{main}, and take $Z=Z(I_n,M_x)$ with $\pi\colon Z\to I_n$.

Note that we above described a class of examples $Z=Z_n$, which
have surjective, ring-like, fully closed maps $\pi:Z_n\to I_n$
whose point-inverses are homeomorphic to pairwise disjoint
non-degenerate subcontinua of $M$, and each of the continua has a
pair of opposite end points. Thus, the statement \ref{snake}(c)
implies that for each $t\in[0,1]$, $(\pi^n_1\pi)^{-1}t$ belongs to
the class of examples $Z_{n-1}$.

By Proposition \ref{comp}, the composition $\pi^n_1\pi$ is ring-like. Using induction on $n$, Lemma \ref{chain} and the
assertions (a--d) of Example \ref{snake}, we infer that $Z$ is a chainable continuum with a pair of opposite end points.
Using induction, Proposition \ref{partit}, and the assertions (c, e) of Example \ref{snake}, we infer that every partition
$L$ in $Z$ has $\ind L\geq n-1$, and hence, $\ind Z\geq n$. By Theorem \ref{indok} applied to $\pi$, $\Indo Z\!=\!n+1$.

Similarly as in the proof of Theorem \ref{kuka}, we shall show
that $Z$ is a weak Cook continuum. Take a continuum $P\subset Z$
and a map $f\colon P\to Z$ with \mbox{$P\cap fP=\emptyset$}.
Suppose that $f$ is not constant. Our first claim is that {\em
$\pi P$ is not a single point and $f$ has a factorization
$f=g(\pi|P)$, where $g\colon\pi P\to Z$ is continuous.} Indeed, if
$Q\subset P$ is a metrizable continuum, then $fQ$ is metrizable,
and is contained in some point-inverse $\pi^{-1}x$, where $x\in
I_n$. Hence, $Q$ and $fQ$ are homeomorphic to disjoint subcontinua
of $M$, and the restriction $f|Q$ must be constant. Thus, $P$ is
not metrizable, $\pi P$ is not a single point by Theorem
\ref{main}(f), and moreover, $f|\pi^{-1}t$ is constant for every
$t\in\pi P$. The map $g\colon\pi P\to Z$, $gt=f\pi^{-1}t$, is well
defined and\linebreak continuous. Our claim ensures that the below
set is non-empty, and we define
\[k_0=\min\{k\colon\mbox{$f$ has a factorization $f=g_k(\pi^n_k\pi |P)$, where $g_k\colon\pi^n_k\pi P\to Z$}\}.\]
Now, observe that, {\em if $Q\subset\pi^n_{k_0}\pi P\subset
I_{k_0}$ is a metrizable continuum, then $g_{k_0}|Q$ is constant.}
Indeed, $Q$ is an arc, $g_{k_0}Q\subset Z$ is metrizable, and is
contained in some point-inverse $\pi^{-1}x$, where $x\in I_n$. As
$M_x$ does not contain arcs, $g_{k_0}|Q$ must be constant. This
shows that $\pi^n_{k_0}\pi P$ is non-metrizable, $k_0>1$, and
$\pi^n_{k_0-1}\pi P$ is not a single point. As $\pi^{k_0}_{k_0-1}$
is atomic, every point-inverse $(\pi^{k_0}_{k_0-1})^{-1}t$,
$t\in\pi^n_{k_0-1}\pi P$, is a terminal continuum, and in
consequence, $(\pi^{k_0}_{k_0-1})^{-1}t\subset\pi^n_{k_0}\pi P$.
By the observation emphasized above, the restriction
$g_{k_0}|(\pi^{k_0}_{k_0-1})^{-1}t$ is constant for every
$t\in\pi^n_{k_0-1}\pi P$. Thus, the map
$g_{k_0-1}\colon\pi^n_{k_0-1}\pi P\to Z$,
$g_{k_0-1}t=g_{k_0}(\pi^{k_0}_{k_0-1})^{-1}t$, is well defined and
continuous. We have $f=g_{k_0-1}(\pi^n_{k_0-1}\pi |P)$, and this
contradicts the definition of~$k_0$. Therefore, $f$ must be a
constant map.

In the same way as in the proof of Theorem \ref{kuka}, one shows
that $Z$ is a Cook continuum.
\end{proof}

\section{Remarks and open problems}

Proposition \ref{partit} or alternatively Theorem \ref{cm} allow us to iterate the constructions in Section \ref{ach} in
order to obtain continua $Z$ with arbitrarily large difference $\ind Z - \dim Z>0$. For example, one can take $Y_x=[0,1]^2$
for every $x\!\in\! [0,1]^2\!=\!Z_1$, apply Theorem \ref{main}, and have $Z_2\!=\!Z([0,1]^2,Y_x)$, $\pi^2_1\colon\! Z_2\!\to
[0,1]^2$, $\dim Z_2=2$. Then, one takes pairwise disjoint non-degenerate subcontinua $C_x$ of Cook's $M_1$ for $x\in Z_2$,
and puts $Z=Z(Z_2,C_x)$ with $\pi\colon Z\to Z_2$. It follows from \ref{comp} and \ref{cm} that $\ind Z\geq 3>2=\dim Z$.
$Z$~is a Cook continuum by the same argument as in the proof of Theorem \ref{cooksnake}.

It follows from Remark 1.8.(2) that, {\em if\/ $\mathsf{CH}$ is
true, then all the examples of continua $Z$ constructed in
Sections \ref{ach}--\ref{secchain} can be perfectly normal,
hereditarily separable, and have\/ $\ind Z\!=\!\Ind Z\!=\!\Indo
Z$.} In Section \ref{secchain}, instead of Chatyrko's continua
$I_n$ one should use perfectly normal chainable continua by A.A.\
Odintsov~\cite{odin}.

We suggest the following open problems.

In most of our examples of continua~$Z$ we have got the annoying
difference between a lower bound of $\ind Z$ and the exact value
of $\Indo Z$.

\begin{statem}\label{ind-domain}
{\bf Question.} Suppose that $f\colon X\to Y$ is a surjective, fully closed ring-like map from a continuum $X$, $Y$ is the
interval\/ $[0,1]$ or another metrizable one-dimension\-al continuum, and every point-inverse of $f$ is a metrizable
one-dimensional\/ {\em [}$n$-dimensional\hspace{.4ex}{\em ]} continuum. Can it happen that\/ $\ind X=1$ {\em
[}\hspace{-.1ex}respectively: $\ind X=n$ or even\/ $\Ind X=n${\em ]}?
\end{statem}

By Theorem \ref{indok}, the above continuum $X$ must have $\Indo
X=2$ [respectively: $\Indo X=n+1$].

\begin{statem}
{\bf Question.} Do there exist\/ {\em [\hspace{-.03ex}}hereditarily indecomposable\hspace{.15ex}{\em ]} Cook~contin\-ua
whose $\dim =n$ for\/ $n\geq 3$ {\em [\hspace{-.02ex}}respectively: $n\geq 2${\em ]}?
\end{statem}

Such continua do not exist in the metric case, see \cite{krz1} and
\cite{mac2}.

\begin{statem}
{\bf Question.} Does there exist a continuum\/ {\em [}\hspace{-.1ex}a dense-in-itself zero-dimensional compact
space\hspace{.15ex}{\em ]} no two of whose disjoint infinite closed subsets are homeomorphic?
\end{statem}

Fe\-dor\-chuk \cite{fed_hered_n-dim1,fed_hered_n-dim2,fed1} and
Chatryko \cite{chat1} constructed examples of hereditarily
$n$-dimensional (with respect to the covering dimension, $n>0$ is
arbitrary) continua. The example in \cite{fed_hered_n-dim2}
($\mathsf{CH}$ assumed) does not have infinite zero-dimensional
closed subspaces, and if it could be Anderson-Choquet, it would
not contain a pair of disjoint homeomorphic infinite closed
subsets.

\end{document}